\documentclass[12pt,fullpage]{article}
\newcommand{\K}{\rm K}
\newcommand{\dg}{{\rm deg}_x}
\newcommand{\tdg}{{\rm trailDeg}_x}
\newcommand{\lc}{{\rm lc}_x}
\newcommand{\pp}{{\rm primpart}}
\newcommand{\flc}{{\rm flc}_x}
\newcommand{\ftc}{{\rm ftc}_x}
\newcommand{\res}{{\rm res}_x}
\newcommand{\cols}{\rm cols}
\newcommand{\rows}{\rm rows}
\newcommand{\tc}{{\rm tc}_x}
\renewcommand{\S}{\rm S}
\newcommand{\genPRem}{\rm genPRem}
\newcommand{\prem}{\rm prem}
\newcommand{\cont}{\rm cont}
\newcommand{\tprem}{\rm tprem}
\newcommand{\relativeSize}{\rm relativeSize}

\begin{document}
\title{Generalized subresultants and generalized subresultant algorithm}
\author{Petr Glotov\\
%Moscow State University\\
 %                      Russia\\
                        pglotov@yahoo.com}
\date{}
\maketitle
\begin{abstract}
In this paper we present the notions of trail (pseudo-)division, 
generalized subresultants and generalized subresultant algorithm.
\end{abstract}

\section{Trail pseudo-division}
We will work in some polynomial ring $\K[x]$.
First, we define {\it full} polynomial as one with non-zero trail coefficient.
All through the paper we will deal with full polynomials. Let we have two
full polynomials $f$ and $g$.
Now we describe the process of trail (pseudo-)division.
The usual (pseudo-)division can be illustrated in the following scheme,
where the coefficients of both polynomials are written from the left to the right by decreasing
powers of $x$ and
the second, third, etc. lines subtracts one after one from the first 
after some multiplications ("$*$" mean some coefficient, no other comments are 
needed there):\\
\begin{tabular}{ccccccccr}
&*&*&*&*&*&*& & $f$\\
&*&*&*&*& & & & $x^2g$\\
&~&*&*&*&*& & & $xg$\\
&~&~&*&*&*&*& & $g$
\end{tabular}\\
We introduce the trail (psudo-)division, which
can be analogously illustrated:\\
\begin{tabular}{ccccccccr}
&*&*&*&*&*&*& & $f$\\
& & &*&*&*&*& & $g$\\
& &*&*&*&*& & & $xg$\\
&*&*&*&*& & & & $x^2g$
\end{tabular}\\
Here the eliminations perform in the trail part of "bigger" polynomial $f$ by
the trail part of the "smaller" one $g$.
In each step one get polynomial from the ideal $(f,g)$ (it has zero in some
lower terms).
After removing of some maximal possible power of $x$
($x \notin (f,g)$, if $\dg(\gcd(f,g))>0$ as $f$ and $g$ are full
and if $\gcd(f,g)=1$ then resulting polynomial is obviously belongs to $(f,g)$)
one can get the polynomial of degree less than $g$ has.

This way of division is usefull in the case of pseudo-division.
In the usual pseudo-division the first polynomial $f$ is multiplied by some degree of the
leading
coefficient of the second one $g$.
In the trail division it is multiplied by some degree
of the trail coefficient of $g$.
If it is "less" in some sense than the leading coefficient
then the resulting trail pseudo-remainder will have "smaller" 
coefficients than the
usual one. 
Analogously in the case of division (not pseudo) the dividing of $g$ can be
performed by "smaller" term.
The $tprem(f,g)$ will denote the trail pseudo remainder of $f$ and $g$.
The remark: in general we can change the "place" of elimination (e.g first
vanish the leading coefficient, then the trail, again the trail, etc.).

If $h$ is full polynomial than $h^*$ will denote the reverted polynomial (e.g.
$(5x^4+4x^3+3x^2+2x+1)^*=x^4+2x^3+3x^2+4x+5)$.
The following formula is valid up to multiplying by some power of $x$ :
$\tprem(f,g)=(\prem(f^*,g^*))^*$. To prove this formula one need to 
place the mirror near the scheme for trail pseudo-remainder. Then in the mirror
one will see the process of finding usual pseudo-remainder of reverted polynomials.

Our next goal is to develop the algorithm analogous to the subresultant algorithm
 for gcd computation \cite{Knuth} with using trail pseudo-remainders.
For this purpose we fix here the generalized algorithm for pseudo-remainders
$genPRem$: this algorithm gets as input two full polynomials $f$ and $g$,
$\dg g \le \dg f$.
In output it produces the polynomial $r$ together with the following six values:
$r$, $\delta$, $\lambda$, $g$, $\bar g$, $w$.
$r$ is the full part of trail or usual pseudo remainder depending on which way
is better
(usaul pseudo-remainder
algorithm doesn't exclude "superflous" powers of $x$); $\delta=\dg f-\dg g$;
$\lambda=\tdg(\prem(f,g))$ 
if usual pseudo-division is used and
$\tdg(\prem(f^*,g^*))$ if trail pseudo-division was performed;
$g$ and $\bar g$ are  $\lc g$ and $\tc g$ or $\tc g$ and $\lc g$ depending 
on the way of division: first pair in the usual case and the second in the trail
one; $w$ is marker of kind of division: $lead$ or $trail$.
Formally we can write it in the following way:\\
~\\
Algorithm genPRem\\
Input: $u,v$ are full polynomials, $\dg(u)\ge\dg(v)$\\ 
Output: the generalized pseudo-remainder\\
~\\
\medskip
if $\relativeSize(\lc(u))\le\relativeSize(\tc(u))$ then\\
\hbox{~~~}$w:=\prem(u,v)$;\\
\hbox{~~~}return($w/x^{\tdg(w)}$,$\dg(u)-\dg(v)$,$\tdg(w)$,$\lc(u)$,$\tc(u)$,$lead$)\\
else\\
\hbox{~~~}$w:=\prem(u^*,v^*)$;\\
\hbox{~~~}return($(w/x^{\tdg(w)})^*$,$\dg(u)-\dg(v)$,$\tdg(w)$,$\tc(u)$,$\lc(u)$,$trail$)\\
fi;\\

Here $\relativeSize$ is a integer characteristic of some term which
says how big it is. For example, the amount of memory which takes the term
can be used.

\section{Generalized subresultant algortihm and generalized subresultants}

Let $f$, $g$ be the initial {\it full} polynomials,
$\dg f\ge \dg g$. Let 
$u_1=\tilde u_1=\bar u_1=f$, $u_2=\tilde u_2=\bar u_2=g$, $\tilde u_3, \tilde u_4,
\dots$ be the sequence of generalized
remainder: $\tilde u_i=\genPRem(\tilde u_{i-2},\tilde u_{i-1})$.
Of course, the elements of this sequence contains removable factors, we need
this sequence just to define the sequence $\delta_i$:
we denote $\delta_i=\dg{\tilde u_{i+1}}-\dg{\tilde u_i}$, 
$\S_m^n=\sum_{i=m}^n\delta_i$,
$\lambda_i=\lambda-value\ of\ \genPRem(\tilde u_{i-1},\tilde u_i)$.
As in the subresultant algorithm we will investigate the determinants of matrices
which consist of coefficients of polynomials $x^\alpha f$, $x^\beta g$: let 
$$
M_k=
\pmatrix{
x^kf\cr
x^{k-1}f\cr
\cdots\cr
f\cr
x^{k+\delta_1}g\cr         
x^{k+\delta_1-1}g\cr
\cdots\cr
g\cr
}=\pmatrix{
*&*&*&*&*&\cdots&*&*\cr
 &*&*&*&*&*&\cdots&*&*&\cr
 & &\ddots\cr
 & & &*&*&*&*&*&\cdots&*&*&\cr
*&*&*&\cdots&*&*\cr
 &*&*&*&\cdots&*&*&\cr
 & &\ddots\cr
 & & & & &*&*&*&\cdots&*&*&\cr
}$$
$$
\qquad\qquad\qquad\qquad\qquad
\underbrace{\qquad\qquad}_{{\rm fixed}\ a\ {\rm cols}}
\qquad\qquad\qquad\qquad
\underbrace{\qquad\qquad\qquad\qquad}_{{\rm fixed}\ 2(k+1)+\delta_1-a-1\ {\cols}}
\quad
$$
where $k<\dg g$.
We denote by $(u_1,u_2)^j$  the polynomial whose coefficients are obtained
by fixing some $a<\rows(M_j)$ columns in the left part of $M_j$, 
$\rows(M_j)-a-1$ columns in the right
side of $M_k$ and calculating the determinant which consist of fixed columns and
one non-fixed column. 
We are not specifying what is $a$ in our consideratuions, we know that there is {\it some} $a$.
$(u_1,u_2)^{S_2^k}$ will be denoted as $\bar u_{k+2}$.
Our goal is to express such polynomial via taking generalized pseudo-remainders.
Most of equations bellow will be true up to the sign -- the sign is
not important in our considerations and it's checking is redundant.
The following relation will be usefull for us: it describes what is happened 
when we perform the generalized pseudo-division in the matrix:

\begin{equation}
(\bar u_1,\bar u_2)^{S_2^k}=
{ {{{\bar g}_2}^{\lambda_2}g_2^{\delta_1+\delta_2-\lambda_2}} \over
{g_2^{(\delta_1+1)(S_2^k+1)}}}(u_2,u_3)^{S_3^k}=
  {{{\bar g}_2}^{\lambda_2}\over
  {g_2^{\lambda_2}}}{1\over{g_2 (g_2^{\delta_1})^{\delta_2}}}
  {(u_2,u_3)^{S_3^k}\over{(g_2^{\delta_1+1})^{S_3^k}}}.\label{transformation}
\end{equation}
From this formula we see, for example, that $\bar u_4=(u_1,u_2)^{S_2^2}=
  {{{\bar g}_2}^{\lambda_2}\over
  {g_2^{\lambda_2}}}{1\over{g_2 (g_2^{\delta_1})^{\delta_2}}}
  {u_4}$, where $u_4=\genPRem(\bar u_2,\bar u_3)$, $\bar u_3=\genPRem(u_1,u_2)$,
i.e. we know what the expression can be removed from $u_4$.
Let $u_i$ denote $\genPRem(\bar u_{i-2},\bar u_{i-1})$.
We want to determine how $u_i$ linked with $\bar u_i$. Let us fix the
number $k$. Then we can write down the following sequence of equations:
$$
\bar u_{k+1}=(\bar u_1,\bar u_2)^{S_2^{k-1}}=G_4^{k+1}(\bar u_2,\bar u_3)^{S_3^{k-1}}=\cdots
$$
$$
G_{i+1}^{k+1}(\bar u_{i-1},\bar u_i)^{S_i^{k-1}}=\cdots=
G_{k+1}^{k+1}(\bar u_{k-1},\bar u_k)^0=G_{k+1}^{k+1}u_{k+1}.
$$

Now we proceed the same transformations with $k$ instead of $k-1$ and simultaneously
we will express $G_l^{k+2}$ via $G_l^{k+1}$ using the~(\ref{transformation}):
\begin{eqnarray}
\bar u_{k+2}=(\bar u_1,\bar u_2)^{S_2^k}=
{{(1/G_3^3)^{\delta_k}}\over{(g_2^{\delta_1+1})^{\delta_k}}}G_4^{k+1}(\bar u_2,\bar u_3)^{S_3^k}=\cdots\cr
={\Biggl({{\prod{1/G_j^j}}\over{\prod{g_{j-1}^{\delta_{j-2}+1}}}}\Biggr)^{\delta_k}}G_{i+1}^{k+1}(\bar u_{i-1},\bar u_i)^{S_i^k}=\cdots
={\Biggl({{\prod{1/G_j^j}}\over{\prod{g_{j-1}^{\delta_{j-2}+1}}}}\biggr)^{\delta_k}}G_{k+1}^{k+1}(\bar u_{k-1},\bar u_k)^{\delta_k}\cr
{=\Biggl({{\prod{1/G_j^j}}\over{\prod{g_{j-1}^{\delta_{j-2}+1}}}}\Biggr)^{\delta_k}}
{{{\bar g}_k}^{\lambda_k}\over{{g_k}^{\lambda_k}}}
{{(1/G_{k+1}^{k+1})^{\delta_k+1}}\over{g_k(g_k^{\delta_{k-1}})^{\delta_k}}}G_{k+1}^{k+1}
(\bar u_k,\bar u_{k+1})^0=G_{k+2}^{k+2}u_{k+2}.\nonumber
\end{eqnarray}
Hence
$$
G_{k+2}={{\bar g}_k^{\lambda_k} \over {g_k^{\lambda_k}}}{1\over{g_k}}
\Biggl({1\over{\prod{G_j^jg_{j-1}^{\delta_{j-2}+1}}G_{k+1}^{k+1}g_k^{\delta_{k-1}}}}\Biggr)^
{\delta_k}.
$$
Let us denote the expression with product as $h_{k+2}$:
\begin{equation}
\label{h}
h_{k+2}={\prod{G_j^jg_{j-1}^{\delta_{j-2}+1}}G_{k+1}^{k+1}g_k^{\delta_{k-1}}}.
\end{equation}
%It is wonder that 
$h_{k+2}$ is "integer" as it is equal to the determinant with "integer" 
entries:
\begin{eqnarray}
({\bar u}_1,{\bar u}_2)^{S_2^{k-1}-1}=G_3^3g_2^{\delta_1+1}G_4^{k+1}({\bar u}_2,{\bar u}_3)^{S_3^{k-1}-1}\cr\bigskip
=G_3^3g_2^{\delta_1+1}G_4^4g_3^{\delta_2+1}G_5^{k+1}({\bar u}_3,{\bar u}_4)^{S_4^{k-1}-1}=\cdots\cr\bigskip
=\prod{G_j^jg_{j-1}^{\delta_{i-2}+1}}G_{k-1}^{k+1}({\bar u}_{k-3},{\bar u}_{k-2})^{S_{k-2}^{k-1}-1}
=\prod{G_j^jg_{j-1}^{\delta_{i-2}+1}}G_k^{k+1}({\bar u}_{k-2},{\bar u}_{k-1})^{\delta_{k-1}-1}\cr
=\prod{G_j^jg_{j-1}^{\delta_{i-2}+1}}G_{k+1}^{k+1}g_k^{\delta_{k-1}-1}u_k,\label{reg}
\end{eqnarray}
and taking the leading or trail coefficient 
%corresponding to the way of the last division 
we get $h_{k+2}$.
%$$
%\prod{G_j^jg_{j-1}^{\delta_{i-2}+1}}G_{k+1}^{k+1}g_k^{\delta_{k-1}}=h_{k+2}.
%$$

We can remark here that from the~(\ref{reg}) it follows that 
$({\bar u}_1,{\bar u}_2)^{S_2^{k-1}-1}\sim {\bar u}_k$ and as we know one of its coefficient,
we can compute it from the ${\bar u}_k$.
Analyzing the view of matrices $M_i$, $S_2^{k-1}< i <S_2^k-1$ (namely, the presence
of zero's on the "leading" or "trailing" diagonals) we see that we can fix the columns in such a way
that $({\bar u}_1,{\bar u}_2)^i=0$ for that $i$, so the structure of the sequence of 
$({\bar u}_1,{\bar u}_2)^i$ is analogue to the one of usual subresultants.

From the~(\ref{h}) it follows the law of $h_k$ transformation:
\begin{eqnarray}
h_{k+2}=h_{k+1}g_{k-1}G_{k+1}^{k+1}g_k^{\delta_{k-1}}=
h_{k+1}g_{k-1}{1\over{h_{k+1}^{\delta_{k-1}}}}{1\over g_{k-1}}
{{{\bar g}_{k-1}^{\lambda_{k-1}}}\over g_{k-1}^{\lambda_{k-1}}}g_k^{\delta_{k-1}}\cr
={{{\bar g}_{k-1}^{\lambda_{k-1}}}\over g_{k-1}^{\lambda_{k-1}}}
{{g_k^{\delta_{k-1}}}\over{h_{k+1}^{\delta_{k-1}-1}}}.\nonumber
\end{eqnarray}

From the considerations above we can derive the algorithms for computing the gcd and
resultants. Bellow we present the algorithm for gcd computation. 
(we present in the style
a la {\bf Algorithm C} from~\cite{Knuth}):\\
~\\
Algorithm C'\\
Input: $f,g$ are polynomials\\ 
Output: the gcd of $f$ and $g$\\
~\\
\medskip
{\bf C'1}. [Reduce to full and primitive.]
%[Initial settings]
 (u,v):=(f,g),
$d:=\gcd(\cont(u),\cont(v))$, $e:=\min(\tdg(u),\tdg(v))$,
replace $(u,v)$ by \\
$(\pp(u)/x^{\tdg(u)},\pp(v)/x^{\tdg(v)})$.
If $\dg(u)<\dg(v)$ then replace $(u,v)$ by $(v,u)$.
Set $h:=1$, $g:=1$, ${\bar g}=1$, $G:=1$, ${\bar G}:=1$.\\
{\bf C'2}. [General pseudo-remainder.] Apply $\genPRem(u,v)$ and assign
$r$,
$\delta$,
$\lambda$,
$g_2$,
${\bar g}_2$,
$w$.
If $r=0$, then return $dx^ev/\cont(v)$.\\ 
{\bf C'3}. [Adjust remainder.] 
$u:=v$;
$v:=(r{\bar G})/(Ggh^{\delta})$
$g:=g_2$;
${\bar g}:={\bar g}_2$;
$h:={\bar G}g^{\delta}/(Gh^{\delta-1})$;
$G:=g^{\lambda}$;
${\bar G}:={\bar g}^\lambda$;
go to {\bf C'2}\\

In the algorithm for computing the resultant of two full polynomials the algorithm is
almost the same, but one need to return the value of $h$.

For the non-full polynomials the following formula for resultant can be used: 
$\res(xu,v)=\tc(v)\res(u,v)$ (up to the sign, of course).

\section{Implementation}
The algorithms for gcd and resultant computing above was implemented with the Axiomxl 
computer algebra system,
which allows to get an efficient executing code. 
As a coefficient ring it was used the ring of polynomials ${\rm Z}[y]$.
In Axiomxl there are two different structures for dense and sparse
polynomials.
As a $\relativeSize$ it was used the degree for dense polynomials
and number of non-zero terms for sparse polynomials. 
The results of testing is
the following: in the case of dense polynomials
the algorithm is not slower than the usual subresultant algorithm;
on some examples it is times faster than the usual subresltant algorithm.

\section{One property of generalized subresultants.}
Here we change the notation and will notate the generalized subresultants as $S_k^*$ to 
underline the analogues with usual subresultants.
$S_k^*$ means that we get the generalized subresultant from the matrix for the usual
subresultant $S_k$.
The well known property of usual subresultants is that theire formal leading
coefficients ({\it principal resultants}) $\flc(S_k)$ allows one to check the
degree of gcd~\cite{Mishra}.
 The generalized subresultants have the same property, namely,
the following lemma can be proved:

\noindent{\bf Lemma.}
{Let $S_k^*$ be the sequence of generalized subresultants of two full polynomials
$A$ and $B$. Then $\dg(\gcd(A,B))=d$ iff $(\flc(S_0^*)\ or\ \ftc(S_0^*))=\cdots=
(\flc(S_{d-1}^*)\ or\ \ftc(S_{d-1}))=0$ and $\flc(S_d^*)\ne 0$ (then also 
$\ftc(S_d^*)\ne 0$ and back); here $\flc$ and $\ftc$ are formal leading and
trailing coefficient, they are some determinants.}\\
The content and proof of the lemma is almost analogous to the corollary 7.7.9 
from~\cite{Mishra}. We just make here some remarks.
Everywhere in the previous to the corollary 7.7.9 lemma's
in~\cite{Mishra} the ${\rm PSC}_i$ appears it can be substituted by
formal leading or trailing coefficient of the generalized subresultants.
The big role in the proof plays the equation $A(x)T_j(x)+B(x)U_j(x)=C_j(x)$,
where there is some conditions on the degrees of $T_j(x)$, $U_j(x)$ and $C_j(x)$ and
which holds when formal leading coefficient of $S_j$ is vanishes.
%The proof of existance of such an equation is made in~\cite{Mishra} by 
%indeterminate coefficients with analyzing some determinants.
In the our case this equation will be of the form
$A(x)T_j(x)+B(x)U_j(x)=x^*C_j(x)$, where $x^*$ means some power of $x$.
%The proof can be done in the similar manner as in~cite{Mishra}.
%but this equation is almost obvious if we remark that $S_j^*$ belongs to the 
%ideal $(A,B)$ and has the degree at most $j$ and if the leading or trailing
%coeficient vanish then we can get the polynomial of degree less then $j$.

%First, we introduce the convinient notion of the 
%{\it length $\lt$} of some polynomial p: $\lt(p)=\dg(p)-\tdg(p)+1$.
%The $\gcd(u,v)$ of two full polynomials have the length $\lt(\gcd(u,v))=
%\dg(gcd(u,v))+1$.
%If $\dg(\gcd(u,v))=0$, then it is equivalent with 
%$\res(u,v)=S_0=S_0^*=\lc(S_0^*)=
%\tc(S_0^*)\ne 0$ and everything is obvious. Let $d>0$.
%Then we proceed by induction on $i$ with the sequence
%$S_0^*,\cdots,S_i^*,\cdots,S_{d-1}^*$.
%For $i=0$ from the well known property of resultant it follows
%that $S_0=S_0^*=0$. Suppose that we prove the equtions of the statement
%for $i=k-1$ and $\lc(S_k^*)=0$ (the case $\tc(S_k^*)=0$ is analogous).
%Then $\dg(\gcd(u,v))>k-1$ and $\dg(S_k^* < k)$
%$S_i^*$ has the degree not more than $i$ and belongs to the ideal $(u,v)$,
%so if $S_i^*\ne 0$, then 
%$\dg(S_i^*)\ge d$ and for $i=0..d-1$ $S_i^*=0$.  
%If $\lc(S_d^*)\ne 0$ then 

\section{Acknowlegment}
I would like to thank M.Bronstein for usefull discussions, help with Axiomxl
and Axiomxl itself and E.V.Zima for providing me with the copy of~\cite{Knuth}.
This paper was written in the year of 2000 I think.

\end{document}